\documentclass{gtart}
\usepackage{graphicx}
\usepackage{indentfirst}
\usepackage{amsmath,amsfonts,amsthm,amssymb}
\usepackage{mathrsfs}
\usepackage{diagrams}

\newtheorem{thm}{Theorem}[section]
\newtheorem{lem}[thm]{Lemma}
\newtheorem{cor}[thm]{Corollary}
\newtheorem{prop}[thm]{Proposition}
\newtheorem{conj}[thm]{Conjecture}

\newtheorem*{theorem}{Theorem}

\theoremstyle{definition}
\newtheorem{defn}[thm]{Definition}
\newtheorem{rem}[thm]{Remark}
\newtheorem{notation}[thm]{Notation}

\begin{document}

\title{Closed 3-braids are nearly fibred}
\authors{Yi Ni}
\address{Department of Mathematics, Princeton University, Princeton,  New Jersey 08544}
\email{yni@math.princeton.edu}

\begin{abstract}
Ozsv\'ath and Szab\'o conjectured that knot Floer homology detects
fibred links. We will verify this conjecture for closed 3-braids,
by classifying fibred closed 3-braids. In particular, given a
nontrivial closed 3-braid, either it is fibred, or it differs from
a fibred link by a half twist. The proof uses Gabai's method of
disk decomposition.
\end{abstract}


\maketitle

\section{Introduction}

Knot Floer homology was introduced by Ozsv\'ath and Szab\'o
\cite{OSz1}, and independently by Rasmussen \cite{Ra}. The Euler
characteristic of knot Floer homology gives rise to the
Alexander-Conway polynomial. Knot Floer homology contains a lot of
information about the knot or link. For example, it detects the
genera of classical links. Namely, we have the following theorem
due to Ozsv\'ath and Szab\'o (\cite{OSz3}, see also \cite{Ni}).

\begin{theorem}
Suppose $L$ is an oriented link in $S^3$. Let $\chi(L)$ be the
maximal Euler characteristic of the Seifert surfaces bounded by
$L$, and $\mathfrak i(L)=\frac{|L|-\chi(L)}2$, where $|L|$ is the
number of components of $L$. Then
$$\mathfrak i(L)=\max\{i|\widehat{HFK}(L,i)\ne0\}.$$
\end{theorem}

We always refer $\widehat{HFK}(L,\mathfrak i(L))$ as the {\it
topmost term} in the knot Floer homology.

We say an oriented link $L$ is {\it fibred}, if the complement of
$L$ fibers over the circle, and $L$ is the boundary of the fiber.
We say the knot Floer homology of a link is {\it monic}, if the
topmost term is isomorphic to $\mathbb Z$. Ozsv\'ath and Szab\'o
proved that if the link is fibred, then the knot Floer homology is
monic (\cite{OSz2}, see also \cite{Ni}). Thus we naturally have
the following conjecture (see \cite{OSz9}):

\begin{conj}\label{fibred}
If a link in $S^3$ has monic knot Floer homology, then it is a
fibred link.
\end{conj}

Not many interesting cases were tested for this conjecture. In the
knot table, there are exactly thirteen 12-crossing non-fibred
knots, each of which has monic Alexander polynomial, and degree of
the Alexander polynomial precisely gives the genus \cite{FK}.
According to some unpublished computations done by Rasmussen, Ni
and Juh\'asz, these knots do not have monic knot Floer homology.
Moreover, the results in \cite{Ni2} give positive theoretical
evidence to the conjecture.

In this paper, we will verify the conjecture for closed 3-braids.
Our geometric result is

\begin{thm}\label{mainthm}
Suppose link $L\in S^3$ is the closure of a 3-braid, then exactly
one of the following 3 cases happens:
\newline i) $L$ is the 3-component trivial link;
\newline ii) $L$ is fibred;
\newline iii) $L$ or its mirror image is the closure of a nondecreasing positive word $P$.
Moreover, either $P$ is a power of one of $a_1,a_2,a_3$, or $P$ is
started with $a_1$ and ended with $a_3$. Hence in the
corresponding braid diagram, after adding a half twist, we get a
fibred link.
\end{thm}

The exact meaning of case iii) will become clear after
Definition~\ref{nondecr}. Our theorem, together with some simple
computations of knot Floer homology, gives the following

\begin{cor}
A closed 3-braid is fibred if and only if it has monic knot Floer
homology.
\end{cor}

\begin{rem}
Closed 3-braids were classified by Birman and Menasco as links
\cite{BM}. It is proved there that a generic closed 3-braid is
represented by a unique conjugacy class of 3-braids.
\end{rem}

\begin{rem}
We are informed by Alexander Stoimenow that the classification of
fibred closed 3-braids has been obtained in \cite{S2}, with the
assistance of Hirasawa and Murasugi.
\end{rem}

The paper is organized as follows: In Section 2, we will compute
the topmost terms in the knot Floer homology of closed 3-braids.
The computation uses a result of Xu \cite{Xu}. In Section 3, we
apply Gabai's method of disk decomposition to prove Theorem
\ref{mainthm}.

\noindent{\bf Acknowledgements.} We wish to thank David Gabai,
Andr\'as Juh\'asz and Zolt\'an Szab\'o for some helpful
conversations. We also wish to thank Joan Birman for some comments
on the paper. We are especially grateful to Jacob Rasmussen, who
pointed out a crucial mistake in an earlier version of this paper,
and to Xingru Zhang, from whose lecture the author learned Xu's
work on 3-braids.

The author is partially supported by the Centennial fellowship of
the Graduate School at Princeton University.

\section{Knot Floer homology of closed 3-braids}

In this section, we will compute the topmost terms in the knot
Floer homology of closed 3-braids. Although the result can be
deduced from our main theorem, the computation here has its own
interest. And the computation motivates our main theorem.

Of course, the computation becomes easier, if we already know what
the genus of a 3-braid is. Fortunately, this problem was solved by
Xu \cite{Xu}. In order to explain her result, we need some
preparation.

\begin{notation}
$B_3$ denotes the group of 3-braids, $\sigma_1$ and $\sigma_2$ are
the standard generators of $B_3$. Instead of the standard
presentation, we use three generators $a_1=\sigma_1$,
$a_2=\sigma_2$, $a_3=\sigma_2\sigma_1\sigma_2^{-1}$, and the
presentation
$$B_3=<a_1,a_2,a_3|a_2a_1=a_3a_2=a_1a_3>.$$
\end{notation}

\begin{center}
\begin{picture}(300,120)
\put(0,0){\scalebox{0.66}{\includegraphics*[86pt,320pt][540pt,
500pt]{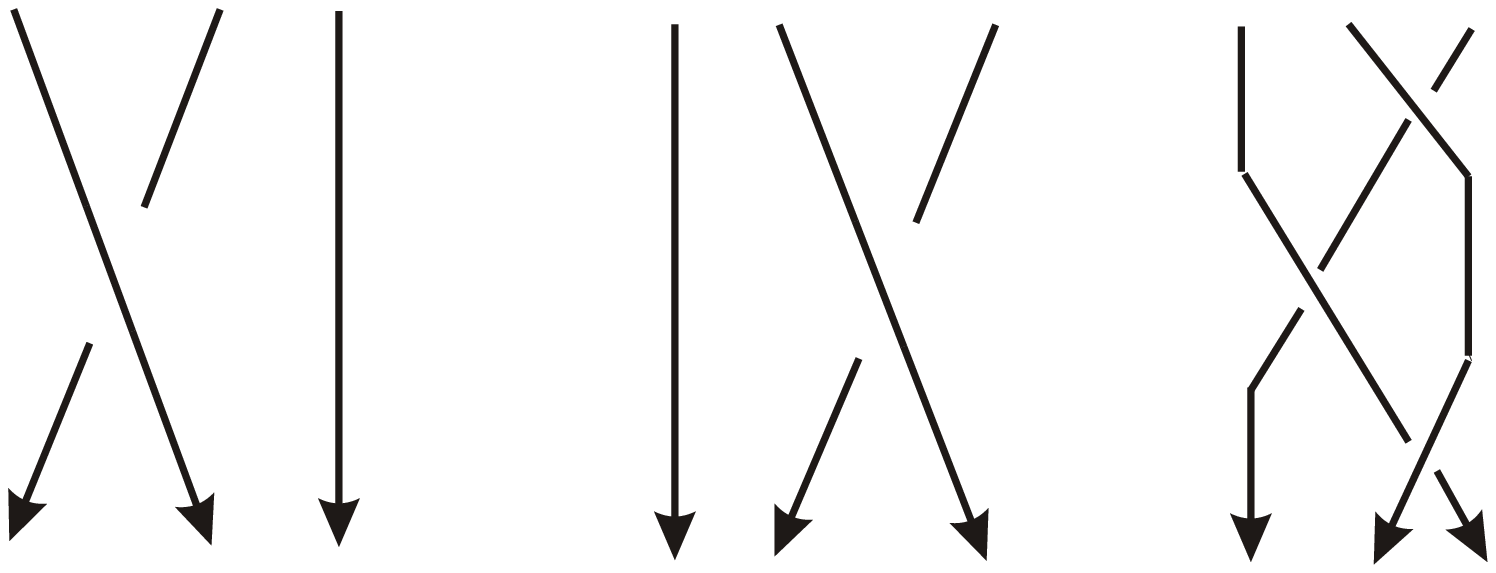}}}

\put(30,1){$a_1=\sigma_1$}

\put(160,1){$a_2=\sigma_2$}

\put(250,1){$a_3=\sigma_2\sigma_1\sigma_2^{-1}$}
\end{picture}

Figure 1 \qquad The generators of $B_3$.
\end{center}

The advantage of such presentation is, one can draw $a_1,a_2,a_3$
cyclically on a cylinder, thus we can permute the roles of
$a_1,a_2,a_3$ cyclically. The reader is encouraged to figure this
out by himself/herself.

If $w\in B_3$, then $\bar w$ denotes its inverse. Let
$\alpha=a_2a_1=a_3a_2=a_1a_3$. The following relations will be
useful to us:
$$a_i\bar a_j=\bar a_{i+1}a_{j+1},\quad \alpha\:\bar a_i=a_{i+1},\quad \bar a_i\alpha=a_{i-1}.$$

\begin{defn}\label{nondecr}
Suppose $P=a_{\varepsilon_1}\cdots a_{\varepsilon_n}$ is a
positive word. We say $P$ is {\it nondecreasing} if for each $j\in
\{1,\dots,n-1\}$, $\varepsilon_{j+1}=\varepsilon_j
\;\text{or}\;\varepsilon_j+1$, where the subscript for $a$ is
understood cyclically. $P$ is {\it strictly increasing}, if for
each $j\in \{1,\dots,n-1\}$, $\varepsilon_{j+1}=\varepsilon_j+1$.
\end{defn}

\begin{thm}[P.J. Xu]\label{Xu}
Every conjugacy class in $B_3$ can be represented by a shortest
word in $a_1,a_2,a_3$ which is unique up to symmetries. The word
has one of the following forms:
\newline
i) $\alpha^kP$;
\newline
ii) $N\alpha^{-k}$;
\newline
iii) $NP$.
\newline
Here $k\ge0$, $\overline N$ and $P$ are nondecreasing positive
words, $P$ or $N$ may be empty.

Moreover, the minimal Seifert surface of the corresponding closed
braid can be constructed from this word.
\end{thm}

We briefly explain how to construct the Seifert surface from a
word $w$. We first resolve the braid to a 3-component trivial
link, bounding 3 disjoint disks. Then for each letter in $w$, one
attaches a twisted band to connect two of the 3 disks. This
surface is called the {\it Bennequin surface} of the word $w$,
denoted by $B_w$. It has Euler characteristic $3-l(w)$, where
$l(w)$ is the length of $w$.

From now on, we also use the word $w$ to denote the corresponding
3-braid, if there is no confusion. Xu's theorem says that, for a
shortest word $w$ as above, $\chi(w)=3-l(w)$.

\begin{rem}\label{reduce}
If the subword $a_ia_i$ appears in $w$, one can replace it by a
single $a_i$ to get a new word $w'$. The Bennequin surface of $w$
is the plumbing of the Bennequin surface of $w'$ with a Hopf band.
It is easy to see the top terms in the knot Floer homology of $w$
and $w'$ are isomorphic as abelian groups. Moreover, the closure
of $w'$ is fibred if and only if the closure of $w$ is fibred
\cite{G1}.

Given a reduced word $w$ in Xu's form, we can apply the previous
``untwisting" operation repeatedly, until we get a word also in
Xu's form, but now the $\overline N$ and/or $P$ are strictly
increasing. We denote this new word by $UT(w)$.
\end{rem}

\begin{thm}\label{FloerHo}
Suppose $L$ is the closure of a 3-braid $w$. $w$ is in the form in
Theorem \ref{Xu}. We consider the word $UT(w)$. If $UT(w)$ is in
the form of $N$ or $P$, and $l(UT(w))=3t+1\;\text{or}\;3t+3$,
($t\ge0$), then
$$\widehat{HFK}(L,\mathfrak i(L))\cong \mathbb Z\oplus\mathbb Z.$$
In other cases, $\widehat{HFK}(L)$ is monic, except when $L$ is
the 3-component trivial link.
\end{thm}

We divide the theorem into several propositions.

\begin{prop}\label{fibredAP}
Suppose $w=\alpha^kP$ is a word in Xu's form, $k>0$. $L$ is the
closure of $w$, then $L$ is fibred with fiber $B_w$. Here $B_w$ is
the Bennequin surface of $w$.
\end{prop}
\begin{proof}
Suppose the first letter in $P$ is $a_1$, then $P=a_1P'$.
$w=\alpha^ka_1P'=\alpha^{k-1}a_2a_1a_1P'$. Hence $B_w$ is the
plumbing of a Hopf band with the Bennequin surface of
$\alpha^{k-1}a_2a_1P'=\alpha^kP'$. By \cite{G1} we can reduce our
problem to $\alpha^kP'$, hence to $\alpha^k$ by induction. Our
conclusion holds since $\alpha^k$ is a torus link.
\end{proof}

\begin{lem}\label{1P}
$L$ is in the form $NP$, $1=l(N)\le l(P)$, then $\widehat{HFK}(L)$
is monic.
\end{lem}
\begin{proof}
We can assume $N=\bar a_2$. We will prove our result by induction
on $l(P)$. When $l(P)=1$, $NP=\bar a_2a_1\;\text{or}\;\bar
a_2a_3$, hence $L$ is the unknot. Now assume $l(P)>1$, and $P$ is
strictly increasing.

If the last letter in $P$ is $a_3$, then $P$ can be written as
$P'a_1a_2a_3$. We have the skein relation for
$$L_-=\bar a_2P'a_1a_2a_3,\quad L_0=\bar a_2P'a_1a_3,\quad L_+=\bar a_2P'a_1\bar a_2a_3.$$
And we have $\bar a_2P'a_1a_3\sim \alpha\bar a_2P'=a_3P'$,
(``$\sim$" denotes conjugacy relation in $B_3$,) hence
$\chi(L_0)\ge\chi(L_-)+3$. In the local picture of the skein
relation, if the two strands in $L_-$ belong to the same
component, then $|L_0|=|L_-|+1$, and $\mathfrak i(L_0)<\mathfrak
i(L_-)$; if the two strands in $L_-$ belong to different
components, then $|L_0|=|L_-|-1$, and $\mathfrak
i(L_0)+1<\mathfrak i(L_-)$. In any case, using the surgery exact
triangle \cite{OSz1}, we get an isomorphism between
$\widehat{HFK}(L_-,\mathfrak i(L_-))$ and
$\widehat{HFK}(L_+,\mathfrak i(L_+))$.

As for $L_+$, we have $\bar a_2P'a_1\bar a_2a_3=\bar
a_2P'a_1a_1\bar a_2\sim \bar a_2^2P'a_1^2$. As we already
mentioned in Remark \ref{reduce}, its knot Floer homology at the
top filtration level is the same as the one of $\bar a_2P'a_1$, to
which we can apply the induction hypothesis.

If the last letter in $P$ is $a_1$, then $P$ can be written as
$P'a_3a_1$. We consider the skein relation for
$$L_-=\bar a_2P'a_3a_1,\quad L_0=\bar a_2P'a_3,\quad L_+=\bar a_2P'a_3\bar a_1.$$
We have $\bar a_2P'a_3\bar a_1\sim P'a_3\bar a_1\bar
a_2=P'a_3\bar{\alpha}=P'\bar a_1$. Length of $P'\bar a_1$ is less
than length of $\bar a_2P'a_3a_1$, hence $\mathfrak
i(L_+)<\mathfrak i(L_-)$. We have $\chi(L_0)=\chi(L_-)+1$. If the
two strands in $L_-$ belong to the same component, then
$|L_0|=|L_-|+1$, and $\mathfrak i(L_0)=\mathfrak i(L_-)$; if the
two strands in $L_-$ belong to different components, then
$|L_0|=|L_-|-1$, and $\mathfrak i(L_0)=\mathfrak i(L_-)-1$. In any
case, we get an isomorphism between $\widehat{HFK}(L_-,\mathfrak
i(L_-))$ and $\widehat{HFK}(L_0,\mathfrak i(L_0))$. Now we apply
the induction hypothesis to $L_0$.
\end{proof}

\begin{prop}\label{NP}
If $L$ is of the type $NP$, $N$ and $P$ are nonempty, then
$\widehat{HFK}(L)$ is monic.
\end{prop}
\begin{proof}
We induct on $l(N)$. The case when $l(N)=1$ is the lemma above.
Now we assume $l(N)\ge2$, we can suppose the first letter in $N$
is $\bar a_3$, then $N=\bar a_3\bar a_2N'$.

If the last letter in $P$ is $a_1$, $P=P'a_1$. Then we consider
the skein relation for
$$L_-=a_3\bar a_2N'P'a_1,\quad L_0=\bar a_2N'P'a_1,\quad L_+=\bar a_3\bar a_2 N'P'a_1.$$
We have $a_3\bar a_2N'P'a_1\sim a_1a_3\bar a_2N'P'=a_3N'P'$. Same
argument as before shows that $\widehat{HFK}(L_0,\mathfrak
i(L_0))\cong\widehat{HFK}(L_+,\mathfrak i(L_+))$. We then apply
the induction hypothesis to $L_0$.

If the last letter in $P$ is $a_2$, $P=P'a_2$. Consider the skein
relation for
$$L_-=a_3\bar a_2N'P'a_2,\quad L_0=\bar a_2N'P'a_2,\quad L_+=\bar a_3\bar a_2 N'P'a_2.$$
$L_0$ can be reduced to $N'P'$, hence we get our conclusion as
before, by applying the induction hypothesis to $L_-=\bar
a_2N'P'a_2a_3$.
\end{proof}

\begin{prop}\label{0P}
If $L$ is in the form $P$, $P=(a_1a_2a_3)^t\;
\text{or}\;(a_1a_2a_3)^ta_1$, $t\ge1$. Then
$$\widehat{HFK}(L,\mathfrak i(L))\cong\mathbb Z\oplus\mathbb Z.$$
\end{prop}
\begin{proof}
Suppose $P=(a_1a_2a_3)^t$, consider the skein relation for
$$L_-=a_2(a_1a_2a_3)^t,\quad L_0=(a_1a_2a_3)^t,\quad L_+=\bar a_2(a_1a_2a_3)^t.$$
$L_-$ can be rewritten as $\alpha P'$, which was considered in
Proposition \ref{fibredAP}, and $L_+$ is of the type considered in
Lemma \ref{1P}. Then $\widehat{HFK}(L_0,\mathfrak i(L_0))$ is fit
into the exact triangle:
\begin{diagram}
\mathbb Z&\rTo &\widehat{HFK}(L_0,\mathfrak i(L_0))\\
\uTo &\ldTo\\
\mathbb Z
\end{diagram}

By \cite{Ni}, $\widehat{HFK}(L_0,\mathfrak i(L_0))\otimes\mathbb
Q$ is nontrivial. One then easily sees that
$$\widehat{HFK}(L_0,\mathfrak i(L_0))\cong \mathbb Z\oplus\mathbb Z.$$

The case when $P=(a_1a_2a_3)^ta_1$ can be reduced to the previous
one by Remark~\ref{reduce}.
\end{proof}

\begin{proof}[Proof of Theorem \ref{FloerHo}]\quad Our theorem now follows
from Remark \ref{reduce}, Proposition \ref{fibredAP}, Proposition
\ref{NP} and Proposition \ref{0P}.\end{proof}

\begin{rem}
With more care, one can get some information of the absolute
grading. For example, in Proposition \ref{fibredAP}, the topmost
term lies at grading level $l(P)+\frac{|L|-1}2$.
\end{rem}

\begin{rem}
Our proof does not really need the fact that the Bennequin surface
of Xu's word is the minimal Seifert surface. This fact can be
proved inductively by our argument.
\end{rem}

\begin{rem}
During the course of this work, we noted the paper \cite{S1}, in
which Stoimenow studied the skein polynomial of closed 3-braids,
also using Xu's theorem. Our result here should be compared with
Stoimenow's work.
\end{rem}

\section{Proof of the main theorem}

\begin{lem}\label{reduce1}
Suppose $w$ is a shortest word for $L$, $w$ is not necessarily in
Xu's form. If the array $\bar a_1a_3a_1a_2$ appears in $w$, then
we can replace the array by $\bar a_1a_2$, thus get a new word
$w'$, with closure $L'$. Then $L$ is fibred with fiber $B_w$, if
and only if $L'$ is fibred with fiber $B_{w'}$.
\end{lem}
\begin{proof}
We draw the local picture of the closed braid near the array $\bar
a_1a_3a_1a_2$ as in Figure~2a. We leave the reader to figure out
the local Bennequin surface. As in \cite{G2}, we thicken $B_w$ to
a sutured manifold $B_w\times I$, and consider its complementary
sutured manifold. In Figure~2b, we draw the suture (as curves) on
the boundary of the handlebody $B_w\times I$. There is an obvious
product disk in the complementary sutured manifold, namely, the
disk bounded by the dashdotted rectangle specified in Figure~2b.
\begin{center}
\begin{picture}(350,193)
\put(0,0){\scalebox{0.7}{\includegraphics*[65pt,380pt][555pt,
650pt]{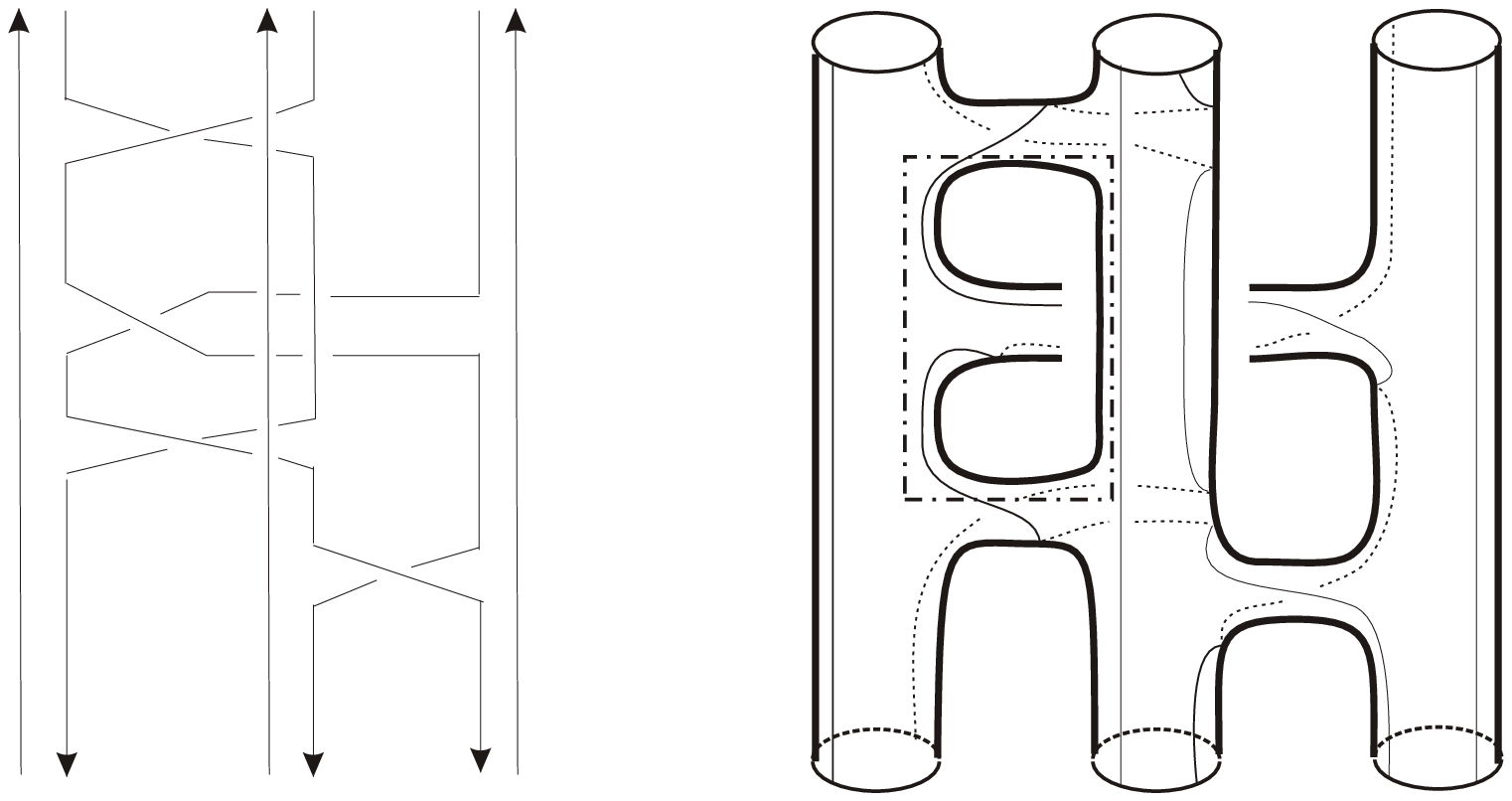}}}

\put(50,5){Figure 2a}

\put(230,5){Figure 2b}
\end{picture}
\end{center}
\begin{center}
\begin{picture}(350,193)
\put(0,0){\scalebox{0.7}{\includegraphics*[50pt,120pt][540pt,
390pt]{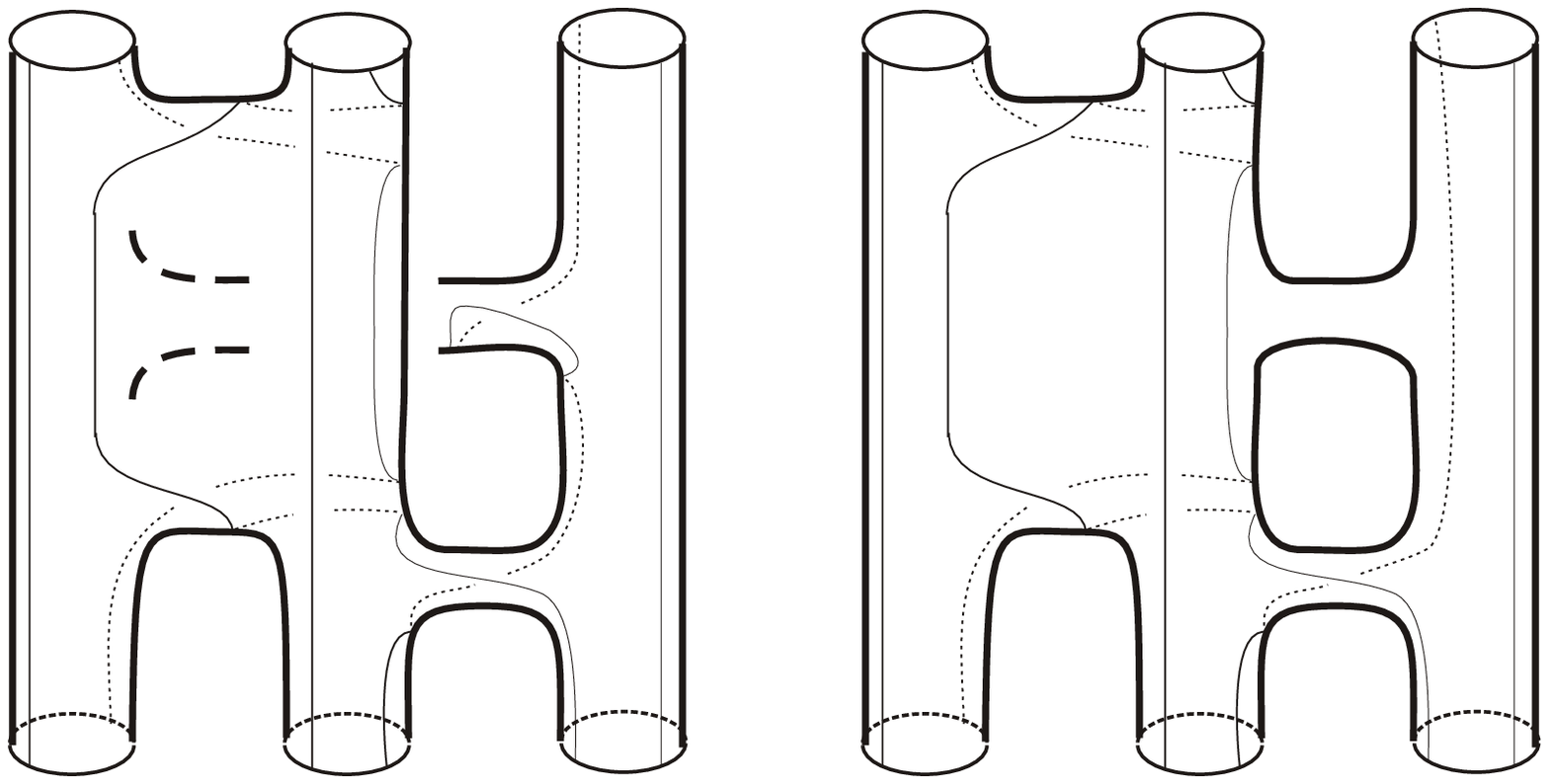}}}

\put(60,1){Figure 3a}

\put(240,1){Figure 3b}

\put(285,80){$D_2$}
\end{picture}
\end{center}
We decompose the complementary sutured manifold along the product
disk, thus get Figure~3a. After an isotopy, we get Figure~3b,
where the product disk $D_2$ is clearer.

Now decompose the complementary sutured manifold in Figure~3b,
thus get Figure~4a. After an isotopy, we get Figure~4b, which is
just the local picture of a Bennequin surface near the array $\bar
a_1a_2$.

Our conclusion holds by Lemma 2.2 in \cite{G2}.
\end{proof}
\begin{center}
\begin{picture}(350,193)
\put(0,0){\scalebox{0.7}{\includegraphics*[50pt,120pt][540pt,
390pt]{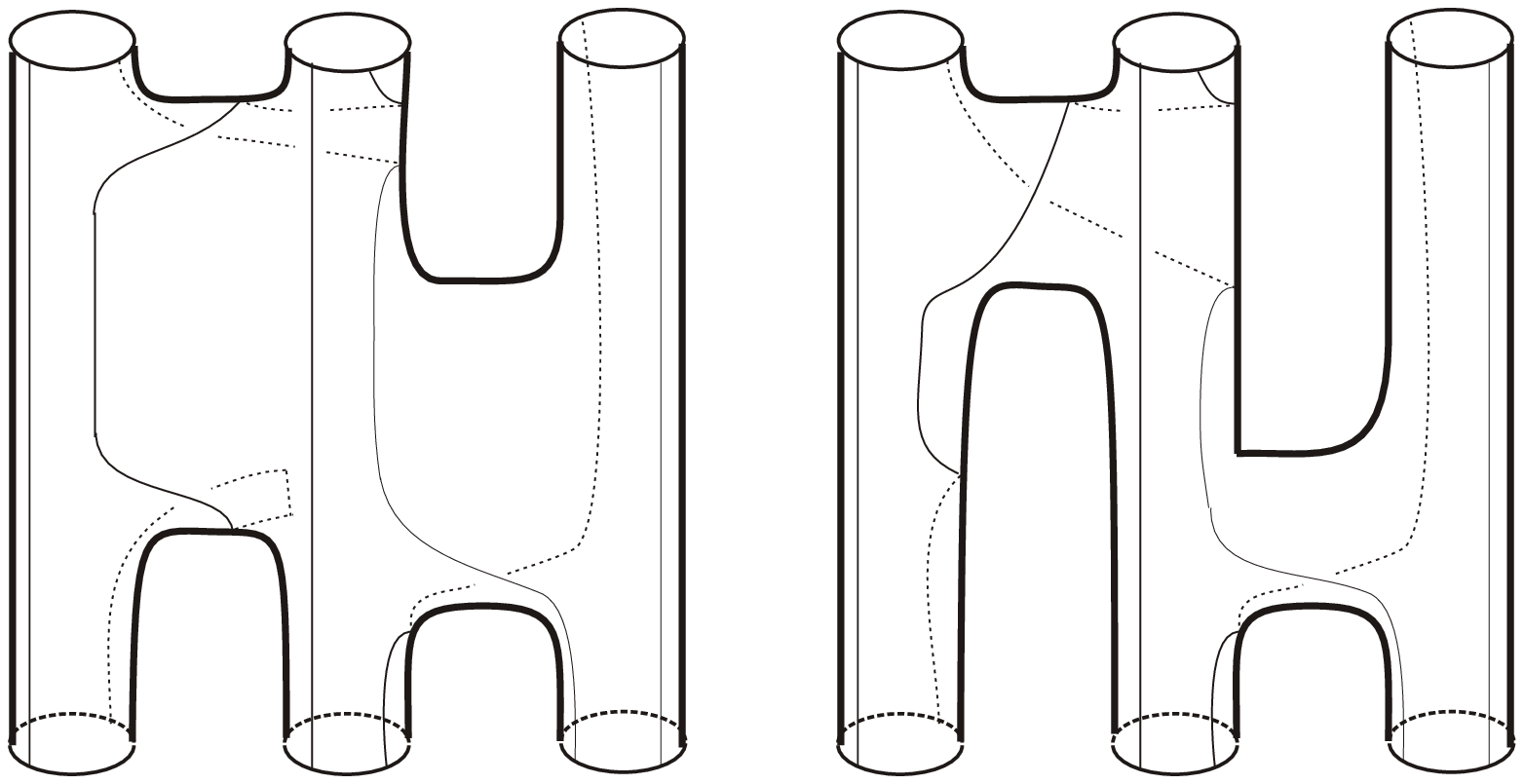}}}

\put(60,1){Figure 4a}

\put(240,1){Figure 4b}
\end{picture}
\end{center}

\begin{lem}\label{reduce2}
Suppose $w$ is a shortest word for $L$, $w$ is not necessarily in
Xu's form. If the array $\bar a_1a_2a_3a_1a_2$ appears in $w$,
then we can replace the array by $\bar a_1a_2$, thus get a new
word $w'$, with closure $L'$. Then $L$ is fibred with fiber $B_w$,
if and only if $L'$ is fibred with fiber $B_{w'}$.
\end{lem}
\begin{proof}
We note that the algebraic relation
$$\cdots\bar
a_1a_2\cdots=\cdots a_3\bar a_1\cdots$$ also gives a local isotopy
of the Bennequin surfaces. We have $\bar a_1a_2a_3a_1a_2=a_3\bar
a_1a_3a_1a_2$. By Lemma \ref{reduce1}, we can replace $a_3\bar
a_1a_3a_1a_2$ by $a_3\bar a_1a_2$. Now $a_3\bar a_1a_2=\bar
a_1a_2^2$, we get our conclusion by Remark \ref{reduce}.
\end{proof}

\begin{prop}\label{fibredNP}
Suppose $w=NP$ is a shortest word in Xu's form for $L$.
$l(N),l(P)>0$. Then $L$ is fibred with fiber $B_w$.
\end{prop}
\begin{proof}
Without loss of generality, can assume $\overline N,P$ are
strictly increasing, and the last letter in $N$ is $\bar a_1$. By
Lemma \ref{reduce1} and Lemma \ref{reduce2}, we can replace $P$ by
one of the following words: $a_3,a_3a_1,a_2,a_2a_3,a_2a_3a_1$.
Then consider $\overline P\:\overline N$. By cyclically permuted
versions of Lemma \ref{reduce1} and Lemma \ref{reduce2}, we can
replace $\overline N$ by a word with length $\le3$. Now there are
only finitely many cases for $NP$ we need to consider. (We note
that $NP$ should be cyclically reduced, this restriction also
reduces the cases.) For these cases, we verify our theorem
directly.
\end{proof}

\begin{proof}[Proof of Theorem \ref{mainthm}]By Proposition \ref{fibredAP}, Proposition \ref{fibredNP}, we only need to
consider the case that $w=P$, $P$ is strictly increasing. The case
that $l(P)\le1$ is easy. If $l(P)=3t+2$, ($t\ge0$,) then it is
conjugated to the form in Proposition \ref{fibredAP}. If
$l(P)=3t+3$ or $l(P)=3t+4$ ($t\ge0$), then one can conjugate the
word so that it is started with $a_1$ and ended with $a_3$. Now
$a_2P$ is fibred by Proposition~\ref{fibredAP}.
\end{proof}

\Addresses

\begin{thebibliography}{H}


\bibitem{BM}{\bf J Birman, W Menasco,} {\it Studying links via closed braids. III. Classifying links which are closed
$3$-braids,} Pacific J. Math. 161 (1993), no. 1, 25--113

\bibitem{FK}{\bf S Friedl, T Kim,}
{\it Thurston norm, fibered manifolds and twisted Alexander
polynomials,}  arXiv:math.GT/0505594

\bibitem{G1}{\bf D Gabai,} {\it The Murasugi sum is a natural geometric
operation,} Low-dimensional topology (San Francisco, Calif.,
1981), 131--143, Contemp. Math., 20, Amer. Math. Soc., Providence,
RI, 1983

\bibitem{G2}{\bf D Gabai,} {\it Detecting fibred links in $S^3$,} Comment. Math. Helv. 61 (1986), no. 4, 519--555

\bibitem{Ni}{\bf Y Ni,} {\it A note on knot Floer homology of links,}
arXiv:math.GT/0506208

\bibitem{Ni2}{\bf Y Ni,} {\it Sutured Heegaard diagram for knots,}
arXiv:math.GT/0507440

\bibitem{OSz1}{\bf P Ozsv\'ath, Z Szab\'o,} {\it Holomorphic disks and knot invariants,} Adv. Math. 186 (2004), no. 1, 58--116

\bibitem{OSz2}{\bf P Ozsv\'ath, Z Szab\'o,} {\it Heegaard Floer homologies and contact
structures,} arXiv:math.SG/0210127

\bibitem{OSz3}{\bf P Ozsv\'ath, Z Szab\'o,} {\it Holomorphic disks and genus bounds,} Geom. Topol. 8 (2004), 311--334 (electronic)

\bibitem{OSz9}{\bf P Ozsv\'ath, Z Szab\'o,} {\it Heegaard diagrams and holomorphic disks,}
Different faces of geometry, 301--348, Int. Math. Ser. (N. Y.),
Kluwer/Plenum, New York, 2004

\bibitem{Ra}{\bf J Rasmussen,} {\it Floer homology and knot complements,} Harvard Thesis, arXiv:math.GT/0306378


\bibitem{S1}{\bf A Stoimenow,} {\it The skein polynomial of closed 3-braids,} J. Reine
Angew. Math. 564 (2003), 167--180

\bibitem{S2}{\bf A Stoimenow,} {\it Properties of closed
3-braids,} preprint

\bibitem{Xu}{\bf P Xu,} {\it The genus of closed $3$-braids,} J. Knot Theory
Ramifications 1 (1992), no. 3, 303--326

\end{thebibliography}
\end{document}